\documentstyle[12pt]{article}
\newlength{\spacing}
\newcommand{\doublespace}{\setlength{\baselineskip}{1.5\spacing}}
\newtheorem{thm}{Theorem}[section]

\newtheorem{prop}[thm]{Proposition}

\def\lam{\lambda }

\def\rar{\to}

\def\al{\alpha}

\def\ep{\epsilon}

\def\lam{\lambda}

\def\cE{{\cal E}}

\def\cN{{\cal N}}

\def\cX{{\cal X}}

\def\today{\ifcase\month\or
  January\or February\or March\or April\or May\or June\or
  July\or August\or September\or October\or November\or December\fi
  \space\number\day, \number\year}

\begin{document}
\begin{titlepage}
\begin{center}
{\bf Number of Edges in Random Intersection Graph on Surface of a Sphere.} \\
\vspace{0.20in} by \\
\vspace{0.2in} {Bhupendra Gupta \footnote{Corresponding Author.
email:gupta.bhupendra@gmail.com, bhupen@iiitdm.in}}\\
Faculty of Engineering and Sciences,\\
Indian Institute of Information Technology (DM)-Jabalpur, India.\\
\vspace{0.1in}
\end{center}
\vspace{0.2in}
\sloppy
\begin{center} {\bf Abstract} \end{center}

\begin{center} \parbox{4.8in}
{In this article, we consider `$N$'spherical caps of area $4\pi p$
were uniformly distributed over the surface of a unit sphere. We
study the random intersection graph $G_N$ constructed by these caps.
We prove that for $p = \frac{c}{N^{\al}},\:c >0$ and $\al >2,$ the
number of edges in graph $G_N$ follow the Poisson distribution. Also
we derive the strong law results for the number of isolated vertices
in $G_N$: for $p = \frac{c}{N^{\al}},\:c >0$ for $\al < 1,$ there is
no isolated vertex in $G_N$ almost surely i.e., there are atleast
$N/2$ edges in $G_N$ and for $\al
>3,$ every vertex in $G_N$ is isolated i.e., there is no edge in
edge set $\cE_N.$}
\\
\vspace{0.4in}
\end{center}

\vspace{0.5in} {\sl AMS 2000 subject classifications}:
\hspace*{0.5in} 05C80, 91D30.
\end{titlepage}
\doublespace
\section{Introduction.}
Random intersection graphs are were introduce in \cite{Singer}, and
defined as:\\
Let us consider a set $V$ with $n$ vertices and another set of
objects $W$ with $m$ objects. Define a bipartite graph $G^*(n,m,p)$
with independent vertex sets $V$ and $W.$ Edges between $v \in V$
and $w \in W$ exists independently with probability $p.$ The random
intersection graph $G(n,m,p)$ derived from $G^*(n,m,p)$ is defined
on the vertex set $V$ with vertices $v_1,v_2 \in V$ are adjacent if
and only if there exists some $w \in W$ such that both $v_1$ and
$v_2$ are adjacent to $w$ in $G^*(n,m,p).$ Also define $W_v$ be a
random subset of $W$ such that each element of $W_v$ is adjacent to
$v \in V.$ Any two vertices $v_1,v_2 \in V$ are adjacent if and only
if $W_{v_1} \cap W_{v_2} \neq \phi,$ and edge set $E(G)$ is define
as
\[E(G) = \{\{v_i,v_j\}: v_i,v_j \in V, \:W_{v_i} \cap W_{v_j} \neq \phi\}.\]
Dudley, \cite{Dudley}, derive the distribution of the degree of a
vertex of random intersection graph. Also show that if $n$ be the
number of vertices and $\lfloor n^{\al}\rfloor$ be the number of
objects, the vertex degree changes sharply between $\al < 1,\;\al
=1$ and $\al>1.$ Bhupendra Gupta \cite{bhupen1} derive the strong
threshold for the connectivity between any two arbitrary vertices of
vertex set $V,$ and determine the almost sure probability bounds for
the vertex degree of a typical vertex of random intersection
graph.\\

{\bf Our Model.}
In this paper we considered the random intersection graph generated
by the spherical caps on the surface of a 3-dimensional
unite sphere.\\
Let $C_1,C_2,\ldots,C_N$ be the spherical caps and
$X_1,X_2,\ldots,X_N$ are their respective centers on the surface of
a unit sphere. Let $X_1,X_2,\ldots,X_N$ are Uniformly distributed
over the surface of unit sphere. Now define a random intersection
graph $G_N$ on the surface of unite sphere, with vertex set $\cX_N =
\{X_1,X_2,\ldots,X_N\}$ and
edge set $\cE_N = \{X_iX_j : C_i\cap C_j \neq \phi,\: i\neq j\}.$\\
The aim of this paper is to investigate the evolution of edges in
the graph $G_N$ with vertex set $\cX_N = \{X_1,X_2,\ldots,X_N\},\:
N=1,2,\ldots,$ where the vertices are independently and uniformly
distributed on the surface of a unit sphere. 
H. Maehara, \cite{maehara} gives the asymptotic results for the
various properties of random intersection graph of random spheriacal
caps on surface of unit sphere. Aslo  Bhupendra Gupta, \cite{bhupen}
gives the strong threshold function $p_0(N) =
o\left(\frac{\log\:N}{N}\right)$ for the coverage of the surface of
a unit sphere by the spherical caps. Bhupendra Gupta shown that for
large $N,$ if $\frac{Np}{\log\:N} > 1/2$ the surface of sphere is
completely covered by the $N$ caps almost surely , and if
$\frac{Np}{\log\:N} \leq 1/2$ a partition of the surface of sphere
is remains uncovered by the $N$ caps almost surely.
\section{Supporting Results.}
Let $C_1,C_2,\ldots,C_N$ be the spherical caps on the surface of a
unit sphere with their centers $X_1,X_2,\ldots,X_N$ and Uniformly
distributed over the surface of unit sphere. We defined a random
intersection graph $G_N$ on the surface of unite sphere, with vertex
set $\cX_N = \{X_1,X_2,\ldots,X_N\}$ and
edge set $\cE_N = \{X_iX_j : C_i\cap C_j \neq \phi,\: i\neq j\}.$\\

Let $p := p(a)$ be the probability that a point `$x$' on the surface
of unit sphere is covered by a specified spherical cap of angular
radius `$a$'. Then the area of the spherical cap of angular radius
`$a$' is equal to $4 \pi p.$\\

\textbf{Poisson Approximation.}\\
Let $\mid \cE \mid$ denote the cardinality of the edge set i.e., the
number of edges in the graph $G_N.$

Define a indicator function
\begin{equation}
\xi_{i} = \left\{
        \begin{array}{ll}
          1, & \hbox{$C_i \cap C_j \neq \phi,\: i \neq j$;} \\
          0, & \hbox{otherwise.}
        \end{array}
      \right.
\end{equation}
that is if $X_i$ is an end point of an edge, then $\xi_i$ is equal
to $1,$ and hence $\mid \cE \mid = \sum_{i\in I}\xi_{i},$ where
$I:=\{i: X_iX_j \in \cE, \: i \neq j\}$ is the index set.

\begin{eqnarray}
E\mid\cE\mid & = & E[\sum_{i=1}^n\xi_i]\nonumber\\
& = & \sum_{i=1}^n E[\xi_i] = {{N}\choose{2}}4p(1-p)\nonumber\\
& = & 2N(N-1)p(1-p) \leq 2 N^2p(1-p).\label{exa_edges}
\end{eqnarray}
%
%
\begin{thm}
(Arratia 1989, \cite{arratia}) Suppose $\xi_i,\: i \in I$ is a
finite collection of Bernoulli random variables.
%
Set $p_i:= E[\xi_i]= P[\xi_i=1],$ and $p_{ij}:= E[\xi_i\xi_j].$ Let
$\lam := \sum_{i\in I}p_{i},$ and suppose $\lam$ is finite. Let
$\mid \cE\mid :=\sum_{i\in I}\xi_i.$ Then
\begin{equation}
d_{TV}(\mid \cE \mid,Po(\lam)) \leq
\min(3,\lam^{-1})\left(\sum_{i\in I}\sum_{j\in \cN_i\setminus
\{i\}}p_{ij}+\sum_{i\in I}\sum_{j\in \cN_i}p_ip_j\right).
\end{equation}
where, $\cN_i$ be the adjacency neighborhood of $i,$ i.e., the set
$\{i\}\cup\{j\in I: X_iX_j \in \cE\}.$\label{arratia}
\end{thm}
\section{Weak Law Results.}
\begin{thm}For $p := p(a) = \frac{c}{N^{\al}},$ where $c>0$ and $\al >2.$ Then sufficiently large $N,$
\begin{equation}
d_{TV}(\mid\cE\mid, Po(\lam)) \rar 0,
\end{equation}
i.e., the number of edges in the graph $G_N$ is a Poisson random
variable with parameter $\lam = \sum_{i\in I}p_{i}<\infty.$
\end{thm}
\textbf{Proof.}
First we consider,
\begin{equation}
p_i  = E[\xi_i] = P[\xi_i=1].\label{edge_pro}
\end{equation}
We know there exists an edge between $X_i$ and $X_j$ iff $C_i\cap
C_j \neq \phi,$ i.e. the distance between $X_i$ and $X_j$ is less
than $2a.$ Now consider another spherical cap $D_i$ centered at
$X_i$ and of radius $2a.$
\begin{eqnarray}
P[\xi_i=1] & = & P[C_i \cap C_j \neq \phi]\nonumber\\
& = & P[X_j \in D_i] = p(2a).\label{e_i}
\end{eqnarray}
Now, from equation (2.1), of Bhupendra \cite{bhupen}, we have
\begin{equation}
p := p(a) = \sin^2(a/2).\label{p}
\end{equation}
Using (\ref{p}) in (\ref{e_i}), we get
\begin{eqnarray}
P[\xi_i =1] & = & \sin^2(a)= \frac{1}{2}(1-\cos(2a))\nonumber\\
& = & 4p(1-p).\label{edge}
\end{eqnarray}
Using (\ref{edge}) in (\ref{edge_pro}), we get
\begin{equation}
p_i  = E[\xi_i] = 4p(1-p).\label{edge_pro_final}
\end{equation}
Now consider
\begin{eqnarray}
p_{ij} & = & E[\xi_i\xi_j]\nonumber\\
& = & 1.P[\xi_i\xi_j = 1] = P[\xi_i=1,\xi_j = 1]\nonumber\\
& = & \sum_{l=1,l\neq i}^{n} \sum_{k=1,k\neq j}^{n}P[(C_i\cap C_l)
\neq \phi,(C_k\cap C_j) \neq \phi]-P[(C_i\cap C_j)
\neq \phi]\nonumber\\
& = & \sum_{l=1,l\neq i}^{n}P[(C_i\cap C_l) \neq \phi]\cdot
\sum_{k=1,k\neq j}^{n}C_k\cap C_j) \neq \phi]-P[(C_i\cap C_j)
\neq \phi]\nonumber\\
& = & \left(4(N-1)p(1-p)\right)^2-4p(1-p)\nonumber\\
& = &
16((N-1)p(1-p))^2\left(1-\frac{1}{4(N-1)^2p(1-p)}\right)\nonumber\\
& \leq & 16((N-1)p(1-p))^2. \label{edge_pro_final1}
\end{eqnarray}
Now by Theorem \ref{arratia}, we have
\[d_{TV}(\mid \cE \mid,Po(\lam)) \leq
\min(3,\lam^{-1})\left(\sum_{i\in I}\sum_{j\in \cN_i\setminus
\{i\}}p_{ij}+\sum_{i\in I}\sum_{j\in \cN_i}p_ip_j\right).\]
Using (\ref{edge_pro_final}) and (\ref{edge_pro_final1}), we get
\begin{eqnarray*}
d_{TV}(\mid \cE \mid,Po(\lam)) & \leq &
\min(3,\lam^{-1})\left(\sum_{i\in I}\sum_{j\in \cN_i\setminus
\{i\}}(4(N-1)p(1-p))^2+\sum_{i\in I}\sum_{j\in
\cN_i}4p(1-p)4p(1-p)\right)\\
& \leq &
\min(3,\lam^{-1})\left(\frac{N(N-1)^3}{2}\left(4p(1-p)\right)^2 +
\frac{N(N-1)}{2}\left(4p(1-p)\right)^2\right).
\end{eqnarray*}
Taking $p = \frac{c}{N^{\al}}$ and $\al >2$ in above, we get
\[d_{TV}(\mid \cE \mid,Po(\lam)) \rar 0, \qquad N\rar \infty.\]
\section{Strong Law Results.}
%
%
\begin{prop}
Let $G_N$ be a random intersection graph. Let $p=
\frac{c}{N^{\al}},$ then
\begin{description}
\item [i.] For $0 < \al < 1,$ there is no isolated vertex in $G_N$ almost
surely.
\item [ii.] For $\al < 2$ at least one isolated vertex in $G_N$ almost
surely.
\item [iii.] For $\al > 3,$ every vertex in $G_N$ is an
isolated vertex.
\end{description}
\label{iso}
\end{prop}
\textbf{Proof.} Let $\cX[B]$ denote that number of vertices of the
finite set point $\cX$ that lies in the set $B.$ Let $D_i$ spherical
cap centered at $X_i$ and of radius $2a.$
\begin{eqnarray*}
P[\mbox{at least one isolated vertex in }G_N] & = &
P[\cup_{i=1}^{N-1}(\cX[D_i]< 1)]\\
& \leq & \sum_{i=1}^{N-1}P[\cX[D_i]< 1]\\
& = & \sum_{i=1}^{N-1}(1-p(2a))^{N-1} = N(1-p(2a))^{N-1}\\
& \leq & (N-1)\exp\left(-(N-1)p(2a)\right)\\
& = & (N-1)\exp\left(-4(N-1)p(1-p)\right),
\end{eqnarray*}
since $p(2a) = 4p(1-p).$ Now taking $p= \frac{c}{N^{\al}},$ we get
\begin{equation}
P[\mbox{at least one isolated vertex in }G_N] \leq
(N-1)\exp\left(-\frac{4(N-1)}{N^{\al}}\left(1-\frac{1}{N^{\al}}\right)\right).
\end{equation}
The above probability is summable for $0 < \al <1,$ i.e.,
\[\sum_{N=1}^{\infty}P[\mbox{at least one isolated vertex in }G_N] < \infty.\]
Then by the Borel-Cantelli's Lemma, we have
\[P[\mbox{no isolated vertex in }G_N,\qquad i.o.] = 1.\]
This implies that for $\al <1$ there is no isolated vertex in $G_N$
almost surely.\\

For the second part of proposition, we consider
\begin{eqnarray*}
P[\mbox{every vertex is an isolated vertex in }G_N] & = &
P[\cap_{i=1}^{N-1}(\cX[D_i]< 1)]\\
& = & \prod_{i=1}^{N-1}P[\cX[D_i]< 1]\\
& = & \prod_{i=1}^{N-1} (1-p(2a))^{N-1} = \left((1-p(2a))^{N-1}\right)^{N-1}\\
& \leq & \left(\exp\left(-(N-1)p(2a)\right)\right)^{N-1}\\
& = & \exp\left(-(N-1)^2p(2a)\right)\\
& = & \exp\left(-4(N-1)^2p(1-p)\right),
\end{eqnarray*}
since $p(2a) = 4p(1-p).$ Now taking $p= \frac{c}{N^{\al}},$ we get
\begin{equation}
P[\mbox{every vertex is an isolated vertex in }G_N] \leq
\exp\left(-\frac{4(N-1)^2}{N^{\al}}\left(1-\frac{1}{N^{\al}}\right)\right).
\end{equation}
The above probability is summable for $\al <2,$ i.e.,
\[\sum_{N=1}^{\infty}P[\mbox{every vertex is an isolated vertex in }G_N] < \infty.\]
Then by the Borel-Cantelli's Lemma, we have
\[P[\mbox{at least one isolated vertex in }G_N,\qquad i.o.] = 1.\]
This implies that for $\al <2$ there is at least one isolated vertex
in $G_N$ almost surely.\\

%
%

%
For the third part of proposition, we consider
\begin{equation}
P[\cE \neq \phi] \leq P[\mid \cE \mid \geq \ep].\label{e1}
\end{equation}
By the Chebyshev's inequality, we have
\begin{eqnarray}
P[\mid \cE \mid \geq \ep] & \leq & \frac{E\mid
\cE\mid}{\ep}\nonumber\\
& \leq & \frac{2}{\ep}N^2p(1-p).
\end{eqnarray}
Taking $p= \frac{c}{N^{\al}},$ we get
\[P[\mid \cE \mid \geq \ep] \leq \frac{2N^{2}}{\ep}\frac{c}{N^{\al}}.\]
Hence from (\ref{e1}), we have
\begin{equation}
P[\cE \neq \phi] \leq \frac{2N^{2}}{\ep}\frac{c}{N^{\al}}.
\end{equation}
The above probability is summable for $\al >3,$ i.e.,
\[\sum_{N=1}^{\infty}P[\cE \neq \phi] < \infty.\]
Then by the Borel-Cantelli's Lemma, we have
\[P[\cE = \phi,\qquad i.o.] = 1.\]
This implies that
\[ \mid \cE \mid = 0,\qquad \mbox{almost surely,}\]
i.e., if $p=\frac{c}{N^{\al}};\:\al>3,$ then there is no edge in the
intersection
graph almost surly, and hence every vertex is an isolated vertex almost surely.\\

\begin{thm}
Let $G_N$ be a random intersection graph. Let $p=
\frac{c}{N^{\al}},$ then for $\al < 1,$ there are at least $N/2$
edges in $G_N$ almost surely. For $\al >3,$ there is no edge in edge
set $\cE_n.$
\end{thm}
\textbf{Proof.} From the Proposition \ref{iso}, we have for $\al <
1,$ there is no isolated vertex in $G_N$ almost surely, i.e., every
vertex is connected with at least one other vertex. This implies
that at least $N/2$ edges in $G_n$ almost surely. \\
For $\al >3,$ every vertex in $G_N$ is isolated almost surely,
implies that there is no edge in $G_N$ almost surely.

\end{document}